\documentclass{amsart}
\usepackage{amsmath,amssymb,amsfonts,color}
\usepackage[all,2cell]{xy}
\SelectTips{cm}{10}

\newtheorem{theorem}{Theorem}[section]
\newtheorem{definition}[theorem]{Definition}
\newtheorem{proposition}[theorem]{Proposition}
\renewcommand{\a}     {{\mathcal A}}
\renewcommand{\c}     {{\mathcal C}}
\renewcommand{\k}     {{\mathcal K}}

\renewcommand{\o}     {{\mathcal O}}
\newcommand{\g}       {{\mathcal G}}
\newcommand{\ba}      {{\mathbf A}}
\newcommand{\bb}      {{\mathbf B}}
\newcommand{\Set}     {{{\sf Set}}}     \let\sets=\Set
\newcommand{\prshv}[1]{{{\sf Set}^{\op{#1}}}}
\newcommand{\ssetsdel}{\prshv{\sdelta}}
\newcommand{\ssets}   {{{\sf Sset}}}    \let\sset=\ssets
\newcommand{\bsets}   {{{\sf B}\ssets}} \let\bset=\bsets
\newcommand{\Cat}     {{\sf Cat}}       \let\cat\Cat

\newcommand{\sdelta}  {{{\Delta}}}

\newcommand{\doscat}  {{2\text{-}\Cat}}

\newcommand{\ldoscat} {{2\text{-}\Cat_{\text{\rm lax}}}}
\newcommand{\xto}[1]  {\xrightarrow{#1}}
\newcommand{\hto}     {\hookrightarrow}
\newcommand{\To}      {\Rightarrow}
\newcommand{\s}[1]    {{\mathbf #1}}
\newcommand{\op}[1]   {{#1^{\text{\rm op}}}}
\newcommand{\LaxFun}  {{\text{\sf LaxFun}}}
\newcommand{\aut}     {{\text{\sf Aut}}}
\newcommand{\Aut}     {{\text{\sf AUT}}}
\newcommand{\ner}     {{\text{\sf Ner}}}
\newcommand{\embed}   {{\text{\bf i}}}

\newcommand{\id}[1]   {{\mathop{\text{\rm id}}\nolimits}_{#1}}
\newcommand{\tdec}    {{\mathop{\text{\rm Tdec}}\nolimits}}
\def\quitar#1{}
\begin{document}

\title{A Full and faithful Nerve for 2-categories}
\author{M.~Bullejos, E.~Faro and V.~Blanco}

%\classno{18G55  (primary),  55S91 (secondary).}

\maketitle

\begin{abstract}
The notion of geometric nerve of a 2-category (Street,~\cite{refstreet}) provides a full and
faithful functor if regarded as defined on the category of 2-categories and  lax 2-functors.
Furthermore,  lax 2-natural transformations between lax 2-functors give rise to homotopies
between the corresponding simplicial maps. These facts allow us to prove a representation theorem
of the general non abelian cohomology of groupoids (classifying non abelian extensions of
groupoids) by means of homotopy classes of simplicial maps.
\end{abstract}

\section{Introduction} \label{intro}

The work presented here began as a search for a definition of nerve of a 2-category which could
be used to prove a representation theorem for the non abelian cohomology of groupoids. This
cohomology, which generalizes the non abelian cohomology of groups given by Dedecker
\cite{Dedecker1964}, is defined in \cite{BlBuFa2004} so as to classify non abelian extensions of
groupoids by ``Schreier invariants'' in a way that generalizes the classical Schreier theory of
groups, (see \cite{Schreier1926a} and \cite{Breen1992}). It was clear that in order to obtain the
desired representation in terms of homotopy clases of simplicial maps, a nerve functor which was
full and faithful was needed. Unfortunately this property was not satisfied by any of the various
notions of nerve of a 2-category existing in the literature. It was only after establishing
Definition \ref{def of geom nerve} that the authors realized that this definition of nerve is
essentially the same as that of nerve of a bicategory of Duskin (see \cite{refduskin}) and that
it agrees on objects with Street's definition of nerve of a 2-category (the so called
\emph{geometric nerve}, \cite{refstreet}).

It is a fairly intuitive fact that, as observed by Grothendieck in the early 1960's, any (small)
category gives rise to a simplicial set (called the \emph{nerve of the category}) and that
functors give rise to simplicial maps between the corresponding nerves. A functor ``nerve'' is
thus obtained from the category $\Cat$ of small categories and functors to that of simplicial
sets, and the fact that this functor is full and faithful has important consequences such as the
fact that $\Cat$ has a model structure whose homotopy category is equivalent to that of
simplicial sets.

Going up one dimension, a similar notion of ``nerve of a 2-category'', associating with each
2-category a simplicial set which captures the essential aspects of the 2-dimensionality of the
structure of the 2-category, becomes useful and necessary in a variety of contexts
(\cite{BuCaFa1998}, \cite{refduskin}, \cite{MoSv1993a}, \cite{refstreet}). In this case, however,
several notions of nerve have been proposed, each with its own applications, but none of them
providing a full and faithful functor.

In this paper we show that if one considers lax 2-functors as the arrows between 2-categories
then, on the one hand, Street's notion of nerve of a 2-category can be given a simple and natural
definition, and, on the other hand, the resulting nerve functor is defined on all lax 2-functors,
that is, it becomes a functor on the category $\ldoscat$ of small 2-categories and lax
2-functors. As a result this nerve functor is full and faithful, allowing us to prove the desired
representation theorem.

The plan of the paper is as follows: Section \ref{gen def of nerve} is an overview of a general
approach to the definition of nerve of the objects of a category. In Section \ref{three nerves}
we review three different notions of nerve of a 2-category and see them as fitting the general
scheme indicated in Section \ref{gen def of nerve}. Special attention is given to the nerve
defined by Street \cite{refstreet} and the fact that it coincides on objects with that of our
Definition \ref{def of geom nerve}. Section \ref{2nerve} is devoted to prove that the geometric
nerve of Definition \ref{def of geom nerve} is full and faithful (Proposition \ref{full 2nerve})
and that to deformations of lax 2-functors there correspond homotopies between their associated
simplicial maps. These results are crucial to establish the representation theorem of non abelian
cohomology of groupoids in Section \ref{repr theor}.

To end this introduction we indicate that all lax 2-functors considered in this paper will be
normal\footnote{Here ``normal'' means ``strictly identity preserving.''}. Thus, the expression
``lax 2-functor'' will be used without further notice to mean ``normal lax 2-functor'', that is,
our convention is that ``a lax 2-functor takes identity 1-cells to identity 1-cells, the
structural natural transformation for identities is an identity, and the structural natural
transformation for composition, $\sigma$, satisfies the additional requirement that for any
1-cell $f:A\to B$,
\begin{equation}
\label{normal}
\sigma_{1_Af} = \sigma_{f1_B} = \id{\bar f},
\end{equation}
(where $\id{\bar f}$ is the identity 2-cell of the image $\bar f$ of $f$).''

\section{General definition of nerve}\label{gen def of nerve}

The classical definition of nerve of a category  reflects the fact that $\sdelta$, the simplicial
category, whose objects are those particular posets that are the finite non-empty linear orders,
$\s1=[0] =\{0\}, \s2=[1] = \{0\leq1\},\dots$ (and whose arrows are all functors or
\emph{monotonic} maps between them), can be regarded as a full subcategory of \Cat. Then, given
that the $n$-simplices of a category $\c$ are just the functors $[n]\to \c$, the fact that the
nerve functor is full and faithful just reflects the fact that $\sdelta$ is an
\emph{adequate}\footnote{The concept of an adequate subcategory was introduced by Isbell in
\cite{Isbell1960}.} subcategory of $\cat$.

By general categorical arguments, given a functor $\embed:\sdelta\to\a$, one can define the
$n$-simplices of an object $A$ of $\a$ as the arrows from $\embed([n])$ to $A$, that is, one can
define
$$
\ner(A)_n = \a\big(\embed([n]),A\big).
$$
The functoriality of $\embed$ provides face and degeneracy operators satisfying the simplicial
identities so that $\ner(A)$ becomes a simplicial set. In this way one obtains a functor
\begin{equation}\label{nerve}
\ner :\a \to \ssetsdel = \ssets
\end{equation}
defined on objects as
\begin{equation}\label{general nerve def}
\ner(A) = \a\big(\embed(-),A\big)
\end{equation}
and on arrows $f:A\to B$ in $\a$ via composition: $\ner(f) = \bar{f}$ is defined as the
simplicial map $$ \ner(A)_n \xto{\ \bar{f}_n\ } \ner(B)_n, \qquad \bar{f}_n(\alpha) =
f\circ\alpha.$$

This construction provides a sensible nerve of the objects of $\a$ if the resulting  functor
\eqref{nerve} is full and faithful (i.e., if $\sdelta$ is adequate in $\a$). If it is not, it may
still be possible to define a sensible nerve of the objects of $\a$ if $\sdelta$ is adequate in
some category $\a^*$ that ``enlarges'' $\a$ in the sense that there is a faithful inclusion
$\a\hto\a^*$. In that case, the nerve of an object $A\in\a$ can be defined as the nerve of $A$ as
an object of $\a^*$. Most of the definitions of nerves that are found in the literature fit this
general scheme. For example, the definitions of nerve of a groupoid and nerve of a 2-groupoid
found in \cite{MoSv1993a} fit this pattern, as they also do all the three definitions of nerve of
a 2-category that we review in Section~\ref{three nerves}.

\section{Three nerve functors for 2-categories}\label{three nerves}

Street in \cite{refstreet} defines a functor $\o:\sdelta \to \omega\text{-}\Cat$ to the category
of (strict) $\omega$-categories and strict functors. Then, he defines the nerve of an
$\omega$-category by means of this embedding, following the general scheme described above. On
the other hand, there is a canonical functor $\omega\text{-}\Cat \to \doscat$ which kills all non
identity cells at dimensions higher than 2 and by composing it with $\o$  an embedding $\bar\o:
\sdelta \to \doscat$ results such that for any simplex $[n]\in \sdelta$ the 2-category
$\bar\o_n=\bar\o([n])$ is the free 2-category in the simplex $[n]$. The resulting nerve functor
$\ner_{\bar\o}:\doscat \to \ssets$ is not full, but the embedding $\bar\o$ has the following
universal property: for any 2-category $\ba$ the set of strict 2-functors $\bar\o_{n} \to \ba$
(or $n$-simplices of $\ba$) is in natural and bijective correspondence with the set of lax
2-functors $[n]\to\ba$. Due to this universal property, the simplicial set $\ner_{\bar\o}(\ba)$
can also be defined in the following way:
\begin{definition}[The geometric nerve]\label{def of geom nerve}
Consider the full and faithful functor
$$
\embed:\sdelta\to\ldoscat
$$
defined by regarding each simplex $[n]$ as a 2-category whose only deformations are the required
identities (this forces  all lax 2-functors $[n]\to[m]$ to be just functors or monotonic maps).
Using this embedding as in \eqref{general nerve def} one gets a nerve functor
$$
\ner : \ldoscat \to \ssets
$$
such that for any 2-category $\ba$, $\ner(\ba)$ is the simplicial set
\begin{equation}
\ner(\ba) = \LaxFun\big(\embed(-),\ba\big), \qquad \vcenter{\xymatrix@C=4mm{ \op\sdelta
\ar[rr]^{\ner(\ba)} \ar[dr]_{\op\embed} & & \sets \\ & \op\ldoscat \ar[ur]_{\LaxFun(-,\ba)} }}
\end{equation}
\end{definition}
By the above universal property, we obtain a commutative diagram
$$
\xymatrix@C-10pt@R+10pt{ \Cat\vphantom{\ldoscat}\ \ar@{^{(}->}[r] \ar[drr]_{\ner} &
\doscat\vphantom{\ldoscat}\ \ar@{^{(}->}[r] \ar[dr]^-(.65){\mskip-10mu\ner_{\bar\o}} & \ldoscat
\ar[d]^-{\ner} \\ & & \sset.}
$$
In other words, the geometric nerve functor just defined on 2-categories agrees with the usual
nerve of categories and with Street's nerve of 2-categories. Furthermore, if  the geometric nerve
is defined for 2-groupoids and that it agrees with the definition of nerve of a 2-groupoid given
by Moerdijk-Svensson in \cite{MoSv1993a}. It should be noted also that the definition of the
geometric nerve agrees also with the definition of nerve of a bicategory given by Duskin
(\cite{refduskin}, p. 238).

Two other useful definitions of nerves of 2-categories can be deduced from an embedding,
$\mathbf{N}:\doscat\hookrightarrow \bset$, of $\doscat$ into the category of bisimplicial sets,
$\bset = \Set^{\op{\sdelta}\times\op{\sdelta}}$ (see diagram \eqref{embed to bisimp} below), and
from two corresponding notions of nerve of a bisimplicial set. The embedding $\mathbf{N}$ of
$\doscat$ into bisimplicial sets  is obtained from the usual nerve of categories by means of
which one can regard  a 2-category as a bisimplicial set in the following way. Given a 2-category
$\ba$ let us write $A_0,A_1$ and $A_2$ for the sets of zero-, one- and two-cells of $\ba$. Then,
we can identify the 2-category $\ba$ with the bisimplicial set $\mathbf{N}(\ba)$ given by the
diagram
\begin{equation}
\label{embed to bisimp}
\vcenter{\xymatrix @R=2pc { & {}\ar@{}[d]|-{\textstyle\vdots} & \ar@{}[d]|-{\textstyle\vdots} &
\ar@{}[d]|-{\textstyle\vdots} \\ \hdots\ \ar@<0.7ex>[r] \ar@<0.2ex>[r] \ar@<-0.3ex>[r]
\ar@<-0.8ex>[r] & {\phantom{AA}} \ar@<0.6ex>[d] \ar[d] \ar@<-0.6ex>[d] \ar@<0.5ex>[r] \ar[r]
\ar@<-0.6ex>[r] & A_2\times_{A_0}A_2 \ar@<0.6ex>[d] \ar[d] \ar@<-0.6ex>[d] \ar@<0.5ex>[r]
\ar@<-0.5ex>[r] & A_1\times_{A_0}A_1 \ar@<0.6ex>[d] \ar[d] \ar@<-0.6ex>[d] \\ \hdots\
\ar@<0.7ex>[r] \ar@<0.2ex>[r] \ar@<-0.3ex>[r] \ar@<-0.8ex>[r] & A_2\times_{A_1}A_2 \ar@<0.5ex>[r]
\ar[r] \ar@<-0.6ex>[r] \ar@<0.5ex>[d] \ar@<-0.5ex>[d] & A_2 \ar@<0.5ex>[r]^{s}
\ar@<-0.5ex>[r]_{t} \ar@<0.5ex>[d]^{s}
\ar@<-0.5ex>[d]_{t} & A_1 \ar@<0.5ex>[d]^{s} \ar@<-0.5ex>[d]_{t} \\
\hdots\ \ar@{=}[r] & A_0\ar@{=}[r] & A_0 \ar@{=}[r] & A_0}}
\end{equation}
whose arrows and columns represent nerves of categories all of which are determined by $\ba$.
This construction gives an embedding $\mathbf{N}:\doscat\hookrightarrow \bset$ which can be
regarded as a sort of ``bi-nerve'' or ``double nerve'' of 2-categories.

This embedding can be composed with two different  ``nerve functors'' of bisimplicial sets. The
first of them is the Artin-Mazur codiagonal (see \cite{artinmazur}) of a bisimplicial set,
$\overline{W}:\bsets\to \sset$, right adjoint to the total dec functor $\tdec: \sset \to \bsets$,
in turn, induced by the ordinal sum functor $+_\text{or}: \sdelta \times \sdelta \to \sdelta$.
Note that $\overline{W}$ can be defined by the general procedure for defining nerves given above
\big(see \eqref{general nerve def}\big), if one uses the embedding of $\sdelta$ into $\bset$
obtained by composing the total dec functor with the Yoneda embedding of $\sdelta$, $\tdec\circ
\mathbf y :\sdelta\hookrightarrow \bset$. The composite $\overline{W}\circ \mathbf{N}$ is a
``nerve'' of 2-categories one of which applications can be seen in in \cite{BuCaFa1998}.

The second ``nerve functor'' of bisimplicial sets comes, again by the general procedure
\eqref{general nerve def}, from the composition $\delta^*\mathbf y$ of the Yoneda embedding of
$\sdelta$ with the left adjoint $\delta^* :\sset\to \bset$ to the functor induced by the diagonal
$\delta:\sdelta\to\sdelta\times\sdelta$, so that we have a further nerve on bisimplicial sets and
on 2-categories $\ner_{\delta}:\doscat\to\sset$. In this case the nerve of a 2-category $\ba$
obtained as the composite $\ner_{\delta}\circ \mathbf{N}$ is the diagonal of the bisimplicial set
$\mathbf{N}(\ba)$.

Each of the three nerves of 2-categories has its own interest in homotopy theory. In
\cite{cal-ceg} and \cite{bull-ceg} it has been proved that the three nerves of a 2-category are
homotopically equivalent. But only the geometric nerve can be extended to lax 2-functors so that
the corresponding nerve on $\ldoscat$ is full and faithful.

\section{The full and faithful geometric nerve functor}
\label{2nerve}

In the following proposition we spell out our definition of the geometric nerve.

\begin{proposition} \label{calculo 2ner en ob} Let $\ba$ be a
2-category, the geometric nerve $\ner(\ba)$ is the following simplicial set: \begin{itemize}
\item Its vertices are the objects of $\ba$, \item 1-simplices are the arrows $A_0\xrightarrow{f}
A_1$ of $\ba$, with faces $$d_0(f)=A_1 \quad \text{and} \quad d_1(f)=A_0,$$ \item 2-simplices are
the diagrams $\Delta = (g,h,f;\alpha)$ of the form $$ 
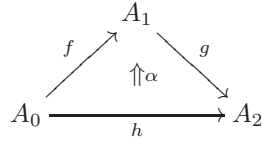
\xymatrix{ & A_1
\ar@{}[d]|-(.55){\hphantom{\alpha}{\textstyle\Uparrow}\alpha} \ar[dr]^{g} \\ A_0 \ar[ur]^{f}
\ar[rr]_{h} & & A_2 } $$ with $\alpha:h\to gf$ a 2-cell in $\ba$, and faces which are the
1-simplices opposite to the indicated vertex, that is,
$$d_0(\Delta)=g, \quad d_1(\Delta)=h,\quad \text{and}\quad
d_2(\Delta)=f,$$ \item 3-simplices are ``commutative'' tetrahedral $\Theta$ of the form \vskip- 2
cm 
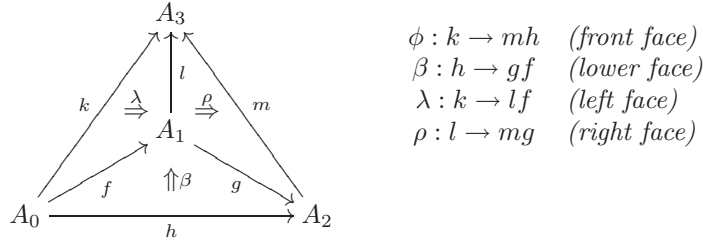
$$ \xymatrix@R-1.5pc@C-0.7pc{ & & A_3 & & \\ & & & & \\ & & & & \\ & & A_1
\ar@{}[dd]|-(.5){\hphantom{\beta}{\textstyle\Uparrow}\beta} \ar[ddrr]_-{g} \ar[uuu]_l
\ar@{}[uuull]|-(.1){\textstyle\stackrel{\lambda}{\To}}
\ar@{}[uuurr]|-(.08){\;\textstyle\stackrel{\rho}{\To}}& & \\ & & & & \\ A_0\ar[uuuuurr]^k
\ar[uurr]_-{f} \ar[rrrr]_{h} & & & &
A_2\ar[uuuuull]_m } \qquad \begin{array}{cl} \\[1.5cm] \phi:k \to mh &
\mbox{(front face)} \\ \beta:h \to gf & \mbox{(lower face)} \\
\lambda:k \to lf & \mbox{(left face)} \\ \rho:l \to mg & \mbox{(right face)} \\ \end{array} $$
where by commutativity of the tetrahedron we mean that the following square of 2-cells commutes
\begin{equation}\label{ctetra} \xymatrix@=5mm{ k \ar[d]_{\lambda}
\ar[r]^-{\phi} & mh \ar[d]^{m\beta} \\ lf \ar[r]^-{\rho f} & mgf }
\end{equation}
The face operators for such tetrahedron are, as in the case of a 2-simplex, the 2-simplices
opposite to the vertex indicated by the operator,

\smallskip

\begin{itemize} \item[] $d_0(\Theta) = \Delta_0$ is
the right face ($\rho$), \item[] $d_1(\Theta) = \Delta_1$ is the front face ($\phi$), \item[]
$d_2(\Theta) = \Delta_2$ is the left face ($\lambda$), and \item[] $d_3(\Theta) = \Delta_3$ is
the lower face ($\beta$). \end{itemize}

\smallskip

\item At higher dimensions $\ner(\ba)$ is coskeletal. \end{itemize}
\end{proposition}

It is also useful to have an explicit description of the simplicial map that the nerve functor
associates to a lax functor.

\begin{proposition} \label{calculo 2ner en ar} Let $F:\ba\to
\bb$ be a lax functor between 2-categories, with structure map $\sigma$. Then $\bar F =
\ner(F):\ner(\ba)\to \ner(\bb)$ is the simplicial map given by: \begin{itemize}

\item[] for 0-simplices $A\in \ner(\ba)_0,\, \bar F_0(A)=F(A)$,

\item[] for 1-simplices $f\in \ner(\ba)_1,\, \bar F_1(f)=F(f)$,

\item[] for 2-simplices $\Delta\in \ner(\ba)_2$, as above, $\bar
F_2(\Delta)$ is obtained from the diagram $$ \xymatrix@C-4mm{ & F(A_1)\ar[ddr]^{F(g)} \\ &
{\hphantom{\textstyle\sigma_{(f,g)}}
{\text{{\large$\Uparrow$}}} {\textstyle\sigma_{(f,g)}}} \\
F(A_0)\ar[uur]^{F(f)} \ar@/^1pc/[rr]^{F(gf)} \ar@/_1pc/[rr]_-{F(h)} & {\scriptstyle
\hphantom{F(\alpha)}{\textstyle\Uparrow} F(\alpha) } & F(A_2) } $$ by the vertical composition of
deformations: $\sigma_{(f,g)}F(\alpha)$, That is $$\bar F_2(g,h,f;\alpha) =
(F(g),F(h),F(f);\sigma_{(f,g)}F(\alpha)).$$

\item[] At higher dimensions $\ner(F)$ is defined in the unique
possible way, using the fact that in dimension 3 and above any simplex is determined by its
faces. \end{itemize} \end{proposition}

Next proposition shows that the construction of the geometric nerve functor is full and faithful.

\begin{proposition}
\label{full 2nerve}
The geometric nerve functor $\ner : \ldoscat \to\sset$ is full and faithful.
\end{proposition}

\begin{proof} We first prove that $\ner$ is full. Let
$\Phi:\ner(\ba)\to\ner(\bb)$ be a simplicial map, then we define a lax functor $F:\ba\to \bb$
acting on objects as the vertex part of $\Phi$: $F(A) = \Phi_0(A)$, on arrows as the edges part
of $\Phi$: $F(f)=\Phi_1(f)$, and on a 2-cell $\alpha:h\to f$ as the 2-cells part of $\Phi$ acting
on the following 2-simplex in $\ner(\ba)$ $$ \xymatrix{{}\ar@{}[d]|-{{\textstyle a =
(1_{A_1},h,f;\alpha)}} \\ {}} \xymatrix{{}\ar@{}[d]|-(.53){{\textstyle\ =}} \\ {}} \xymatrix{ &
A_1 \ar@{}[d]|-(.56){\hphantom{\alpha} {\textstyle\Uparrow}\alpha} \ar@{=}[dr] \\ A_0 \ar[ur]^{f}
\ar[rr]_{h} & & A_1 } $$ $F(\alpha) = \Phi_2(a)$. To complete the definition of $F$ we define its
structure map: given a composable pair of arrows $A_0 \xto{f} A_1 \xto{g}A_2$ in $\ba$, we
consider the 2-simplex in $\ner(\ba)$ given by the diagram
$$ \xymatrix{{}\ar@{}[d]|-{{\textstyle b = (g,gf,f;\id{gf})}} \\ {}}
\xymatrix{{}\ar@{}[d]|-(.57){{\textstyle\ =}} \\ {}} \xymatrix{ & A_1
\ar@{}[d]|-(.56){\!\hphantom{\id{gf}} {\textstyle\Uparrow} {\id{gf}}} \ar[dr]^g \\ A_0
\ar[ur]^{f} \ar[rr]_{gf} & & A_1 } $$ then the structure map $\sigma$ is defined by saying that
$\sigma^{A_0 A_1 A_2}_{(f,g)}$ is the 2-cell that is the interior of the 2-simplex $\Phi_2(b)$.

To see that the nerve functor is faithful we first observe that for any lax functor $F:\ba\to
\bb$ and any 2-cell $\alpha:h\to f$ in $\ba$, if $a = (1_{A_1},h,f;\alpha)$ is the 2-simplex in
$\ner(\ba)$ indicated above, then $\ner(F)_2(a)$ is the 2-simplex given by the diagram
$$
\xymatrix@C-4mm{ & F(A_1) \ar@{}[d]|-(.56){\hphantom{F(\alpha)}{\textstyle\Uparrow} F(\alpha)}
\ar@{=}[dr] \\ F(A_0) \ar[ur]^{F(f)} \ar[rr]_{F(h)} & & F(A_1) }
$$
We also observe that, for any composable pair $(f,g)$ in $\ba$, the structure morphism
$\sigma^{A_0 A_1 A_2}_{(f,g)}$ is the 2-cell that is the interior of the 2-simplex
$\ner(F)_2(b)$, where $b$ is as above. Therefore two lax 2-functors which have the same nerve
coincide on objects, arrows and 2-cells and also have the same structure map.
\end{proof}

The nerve construction not only relates lax 2-functors with simplicial maps but it also relates
lax 2-natural transformations with homotopies. In the next proposition, whose proof we leave to
the reader, we state that lax 2-natural transformations between lax 2-functors give rise to
homotopies between their nerves and that these homotopies satisfy a special property. Remember
that for $n\geq3$ the $n$-simplices of a 2-category are completely determined by their boundary
(the faces) and therefore a homotopy to the nerve of a 2-category is completely determined by its
2-trucation.

\begin{proposition} \label{pseudo nat y homotopies} Let $F,G :\ba \to
\bb$ be two lax 2-functors between 2-categories. Every lax 2-natural transformation $\alpha :
F\to G$ induces a homotopy $h= \ner(\alpha):\ner(G) \to \ner(F)$ such that the image by
$h_1^1:\ner_1(\ba)\to\ner_2(\bb)$ of any 1-simplex in $\ner(\ba)$ is a 2-simplex which is
completely determined by its boundary.
\end{proposition}

\begin{proposition}
\label{gorupoid homotopies}
Let $\ba$ be a 2-category and $\bb$ a 2-groupoid. If $F,G :\ba \to \bb$ are two lax 2-functors,
there exists a homotopy from the simplicial map $\bar F = \ner(F)$ to $\bar G = \ner(G)$ if and
only if there is a lax 2-natural transformation $\alpha: G\to F$.
\end{proposition}

\begin{proof}
One implication is a consequence of Proposition \ref{pseudo nat y homotopies}. To prove the
converse, let $h:\bar F\to\bar G$ be a homotopy. Then, for any object $A$ in $\ba$, $h_0^0(A) :
G(A) \to F(A)$ since $d_0(h_0^0 ) = F$ and $d_1(h_0^0 ) = G$. We define a lax 2-natural
transformation $\alpha$ by defining its $A$-component as $\alpha_A = h_0^0(A)$. In order to
define the structure morphism, let now $f:A\to B$ be a 1-cell in $\ba$. The homotopy identities
imply that we have the following 2-cells of $\bb$: $h_0^1(f):d\to F(f)\alpha_A$ and
$h_1^1(f):d\to \alpha_BG(f)$ with $d = d_1h_0^1(A) = d_1h_1^1(A) : G(A) \to F(B)$
$$
\xymatrix{ & F(A) \ar@{}[d]|-(.55){\hphantom{h_0^1(f)}{\textstyle\Uparrow}h_0^1(f)}
\ar[dr]^{F(f)}
\\ G(A) \ar[ur]^{\alpha_A} \ar[rr]_{d} & & F(B) } \qquad \xymatrix{ & G(B)
\ar@{}[d]|-(.55){\hphantom{h_1^1(f)}{\textstyle\Uparrow}h_1^1(f)} \ar[dr]^{\alpha_B} \\ G(A)
\ar[ur]^{G(f)} \ar[rr]_{d} & & F(B). }
$$
Since $\bb$ is a groupoid, the 2-cell $h_1^1(f)$ is invertible and we can define the structure
map $s(f) = h_0^1(f)\big(h_1^1(f)\big)^{-1}$. The homotopy identities in the next dimension
guarantee the coherence condition so that $(\alpha,s)$ is a lax 2-natural transformation from $G$
to $F$.
\end{proof}

\section{Representation theorem for non abelian
cohomology of groupoids}\label{repr theor}

As an immediate application of the full and faithful geometric nerve of 2-catego\-ries, we  give
a representation theorem for the 2-dimensional non abelian cohomology of groupoids defined in
\cite{BlBuFa2004}, in terms of homotopy classes of simplicial maps.

Given two groups, $G$ and $K$, the (``non abelian'') extensions of $G$ with  kernel $K$ were
classified by Dedecker by means of a two dimensional non abelian cohomology,
$H^2_{\aut(K)}(G,K)$. This is a set of cohomology classes of 2-cocycles defined by taking
coefficients in the crossed module $K\to\aut(K)$ of automorphisms of $K$.\quitar{ Dedecker
\cite{Dedecker1964} goes further to show that it is possible to give a general definition of the
2-dimensional cohomology of a group with coefficients in an arbitrary crossed module.} It is
proved in \cite{BlBuFa2004} that, by regarding the group(oid) $G$ as a (2-discrete) 2-groupoid,
and by taking into account the equivalence between the category of crossed modules and that of
2-groupoids, this cohomology can be represented in terms of lax 2-functors and in fact one has
$$H^2_{\aut(K)}(G,K) =
\LaxFun[G,\Aut(K)]$$ where $\Aut(K)$ is the 2-groupoid associated to the crossed module
$K\to\aut(K)$ and $\LaxFun[G,\Aut(K)]$ denotes the set of connected components of the category
$\LaxFun\big(G,\Aut(K)\big)$ of lax 2-functors and lax 2-natural transformations from $G$ to
$\Aut(K)$. This identification of 2-cocycles with lax 2-functors follows immediately from the
observation  that Dedecker's definition of cocycle and cocycle condition are exactly the
definition of a lax 2-functor and coherence of a lax 2-functor. On this basis \cite{BlBuFa2004}
defines the 2-cohomology of a groupoid $\g$ with coefficients in a family of groups
$\k=\{K_x\}_{x\in\g}$ indexed by the objects of $\g$ as the set
$$H^2_{\aut(\k)}(\g,\k) = \LaxFun_*[\g,\Aut(\k)]$$ of connected
components of lax 2-functors which are the identity on objects. This cohomology, which obviously
generalizes Dedecker's, can also be given in terms of cohomology classes of cocycles and it
classifies general non abelian extensions of a groupoid $\g$ by the family of groups $\k$ indexed
by the objects of $\g$. On the other hand, Propositions \ref{full 2nerve} and \ref{gorupoid
homotopies} above give, as an immediate corollary, the following representation theorem of non
abelian cohomology in terms of homotopy classes of simplicial maps:

\begin{theorem}
\label{repr th}
Given a groupoid $\g$ and any family $\k$ of groups indexed by the objects of $\g$, there is a
bijection
$$
H^2_{\aut(\k)}(\g,\k)\cong \big[\ner(\g),\ner\big(\Aut(\k)\big)\big],
$$
where $\Aut(\k)$ is the 2-groupoid with objects the groups $K_x\in\k$, 1-cells their group
isomorphisms (regarded as functors), and 2-cells the natural transformations between them, and
the square brackets mean homotopy classes of simplicial maps.
\end{theorem}

This theorem can be used, for example, to classify homotopy classes of continuous maps from a
1-type to a 2-type in terms of the said non abelian cohomology of groupoids.

\section{Acknowledgments}

We wish to thank J. Duskin for his valuable comments. Our correspondence with him shows that most
of our results were known to him and are implicitly contained in his work on nerves of
bicategories. Even so we felt that our exposition may be useful for a large audience. We have
used Duskin's clear geometric notation for the $n$-simplices and the corresponding convention for
naming their faces.

\bigskip
\begin{tabular}{lccl}
  M. Bullejos, V. Blanco& && E. Faro \\
  Department of Algebra &&& Department of Appl. Mathematics \\
  University of Granada && &University of Vigo \\
  18071 Granada, Spain &&& 36207 Vigo, Spain \\
  bullejos@ugr.es & &&efaro@dma.uvigo.es \\
  vblanco@ugr.es &&&
\end{tabular}


\begin{thebibliography}{10}

\bibitem{Breen1992}
L.~Breen.
\newblock Th\'eorie de {S}chereier sup\'erieure.
\newblock {\em Ann. Sci. \'Ecole Norm. Sup. (4)}, 25, 1992.

\bibitem{BlBuFa2004}
M.~Bullejos, E.~Faro, and V.~Blanco.
\newblock {\em Categorical non abelian cohomology and the Schreier
theory of groupoids}.  {\tt http://www.ugr.es/\~{}bullejos/schreier.pdf}.

\bibitem{BuCaFa1998}
M.~Bullejos, J.~G. Cabello, and E.~Faro.
\newblock On the equivariant $2$-type of a ${G}$-space.
\newblock {\em J. Pure Appl. Algebra}, 129(3):215--245, 1998.

\bibitem{cal-ceg} M.~Bullejos and A.~M. Cegarra. {\em On the geometry
of 2-categories and their classifying spaces (extended version)}, \newblock {\tt
http://www.ugr.es/\~{}bullejos/geometryampl.pdf}, 2003.

\bibitem{bull-ceg} M.~Bullejos and A.~M. Cegarra. \newblock On the
geometry of 2-categories and their classifying spaces. \newblock {\em K-Theory}, 29:211--229,
2003.

\bibitem{artinmazur} J.~M. Cordier. \newblock Sur les limites
homotopiques de diagrammes homotopiquement coh\'erents. \newblock {\em Compositio Math.},
62:367--388, 1987.

\bibitem{Dedecker1964}
P.~Dedecker.
\newblock Les foncteurs Ext$_{II}$, H$^2_{II}$ et H$^2_{II}$ non
abeliens.
\newblock {\em C. R. Acad. Sc. Paris}, 258 Groupe 1:4891--4895, 1964.

\bibitem{refduskin} John~W. Duskin. \newblock Simplicial matrices and
the nerves of weak $n$-categories i: Nerves of bicategories. \newblock {\em Theory and
Applications of Categories}, 9(10):198--308, 2002.

\bibitem{Isbell1960} J.~R. Isbell. \newblock Adequate subcategories.
\newblock {\em Illinois J. Math.}, 4:541--552, 1960.

\bibitem{Isbell1964} J.~R. Isbell. \newblock Subobjects, adequacy,
completeness and categories of algebras. \newblock {\em Razprawy Mat.}, 36:1--32, 1964.

\bibitem{JoSt1991} A.~Joyal and R.~Street. \newblock Braided tensor
categories. \newblock {\em Advances in Math}, 1(82), 1991.

\bibitem{KeSt1974} G.~M. Kelly and R.~Street. \newblock {\em Review of
the elements of 2-categories}, volume 420 of {\em Lecture Notes in Math.} \newblock
Springer-Verlag, 1974.

\bibitem{Lack2002}
Stephen Lack.
\newblock A quillen model structure for 2-categories.
\newblock {\em K-Theory}, 26(2):171--205, Jun 2002.

\bibitem{MoSv1993a} I.~Moerdijk and J.~A. Svensson. \newblock The
equivariant serre spectral sequence. \newblock {\em Proc. Amer. Math. Soc.}, 118(1):263--277,
1993.

\bibitem{MoSv1993b} Ieke Moerdijk and Jan-Alve Svensson. \newblock
Algebraic classification of equivariant homotopy 2-types, part i.
\newblock {\em J. of Pure and Applied Algebra}, 89:187--216, 1993.

\bibitem{Schreier1926a}
O.~Schreier.
\newblock Uber die erweiterung von gruppen {I}.
\newblock {\em Monatshefte f\"ur Mathematik und Physik}, 34:165--180, 1926.

\bibitem{Schreier1926b}
O.~Schreier.
\newblock Uber die erweiterung von gruppen {II}.
\newblock {\em Abh. Math. Sem. Hamburg}, 4:321--346, 1926.

\bibitem{refstreet} Ross Street. \newblock The algebra of oriented
simplexes. \newblock {\em J. of Pure and Applied Algebra}, 49:283--335, 1987.

\bibitem{refstreetHandbook} Ross Street. \newblock {\em Categorical
structures, Handbook of Algebra Volume 1}, pages 529--577. \newblock Elsevier Science, Amsterdam,
1996.

\end{thebibliography}
\end{document}